\magnification 1200
        \def\R{{\rm I\kern-0.2em R\kern0.2em \kern-0.2em}}
        \def\N{{\rm I\kern-0.2em N\kern0.2em \kern-0.2em}}
        \def\P{{\rm I\kern-0.2em P\kern0.2em \kern-0.2em}}
        \def\B{{\rm I\kern-0.2em B\kern0.2em \kern-0.2em}}
        \def\Z{{\rm I\kern-0.2em Z\kern0.2em \kern-0.2em}}
        \def\C{{\bf \rm C}\kern-.4em {\vrule height1.4ex width.08em depth-.04ex}\;}
        \def\B{{\bf \rm B}\kern-.4em {\vrule height1.4ex width.08em depth-.04ex}\;}
        
        \def\D{{\Delta}}

        \def\z{{\zeta}}

        \def\cW{{\cal W}}

        \def\tau{{\sigma}}

        \
        \vskip 20mm
        \centerline {\bf BOUNDARY CONTINUITY OF COMPLETE PROPER}
        \centerline{\bf HOLOMORPHIC MAPS}
        \vskip 4mm
        \centerline{Josip Globevnik}
        \vskip 4mm
        
	        \noindent \bf Abstract\ \ \rm We show that there is no complete proper holomorphic map 
	        from the disc $\D$ to the bidisc $\D^2$ which extends continuously through $\overline\D$. 
        \vskip 4mm
        Let $\D$ be the open unit disc in $\C$ and let $N\geq 2$. A holomorphic 
        immersion $\varphi \colon\D
        \rightarrow \C^N$ is called {\it complete} if the pullback of the 
        Euclidean metric $\varphi^\ast g$ is a complete metric on $\D$. This is equivalent to saying that
        for every path 
        $p\colon[0,1)\rightarrow \D$ such that 
        $p(t)\rightarrow b\D$ as $t\rightarrow 1$, the composition $\varphi\circ p$ 
        has 
        infinite Euclidean length. Notice that this equivalent statement makes sense as 
        the definition of completeness for 
        general holomorphic maps. 
        
        Answering a question of P.\ Yang [Y], P.\ W.\ Jones [J] was the first to show that 
        there are bounded complete holomorphic 
        immersions. It is now known that given a convex domain $D\subset \C^2$ there is a complete 
        proper holomorphic immersion $\varphi\colon\D\rightarrow D$ [AL]. In [AF] this was
        generalized to the 
        case where $\D$ is replaced with a bordered Riemann surface.
        
        It is a natural question whether, in the case of bounded $D$, there is such a  
        $\varphi$ which extends  
        continuously through $\overline\D$. For instance, if $D$ is a ball, does there exist a complete, 
        proper holomorphic map $\varphi \colon\D\rightarrow D$ which extends continuously through 
        $\overline\D$? This is an open question. In the present note we show that for general bounded convex 
        domains $D$ the answer is no:
        \vskip 2mm
        \noindent\bf PROPOSITION \it\  Let $F=(f,g)\colon\D\rightarrow \D^2$ be a complete, 
        proper holomorphic map. There is no 
        arc $\Lambda\subset b\D$ such that $z\mapsto |f(z)| $ extends continuously to $\D\cup\Lambda$. \rm
        \vskip 2mm
        \noindent\bf Proof.\rm\ \  Suppose that $\Lambda\subset b\D$ is
        an open arc 
        such that $z\mapsto |f(z)|$ extends continuously to $\D\cup\Lambda $. Suppose 
        first that the continuous extension of $z\mapsto |f(z)|$ is identically equal 
        to 1 on $\Lambda$. By the Schwarz reflection principle [R, p.\ 237; p.\ 293, Ex.2]
         $f$ extends 
        holomorphically across $\Lambda$. Thus, for every $e^{i\theta} \in\Lambda$ we have
        $$
        \int_0^1|f^\prime (te^{i\theta})|dt <\infty  .
        \eqno (1)
        $$
        Since the function $g$ is bounded, a result of J.\ Bourgain [B] implies that there is an $e^{i\alpha}
        \in\Lambda $ such that 
        $$
        \int_0^1|g^\prime (te^{i\alpha})|dt <\infty .
        $$
        By (1) it follows that 
        $$
        \int_0^1\sqrt{|f^\prime(te^{i\alpha})|^2 + |g^\prime(te^{i\alpha})|^2}dt<\infty
        $$
        which contradicts the completeness of $F$.
        
        Thus, there is a point $w\in\Lambda$ such that the the extension of 
        $z\mapsto |f(z)|$ at $w$ is less than one. By the continuous
        extendibility of $z\mapsto |f(z)|$ to $\D\cup\Lambda $ there are a 
        closed arc $A\subset \Lambda$, a neighbourhood $\cW $ of $A $ in $\C$ and an $\eta >0$ 
        such that
        $$
        |f(z)|< 1-\eta \ \ (z\in\cW\cap\D).
        $$ 
        Since the map $F$ is proper it follows that if $z\in A$ and $\z\in\D,\ \z\rightarrow z$ then $|g(\z )|
        \rightarrow 1 .$ This means that $z\mapsto|g(z)|$ extends continuously to $\D\cup A$ and that 
        the extension is identically equal to $1$ on $A$. This is impossible by the first part of the proof 
        with the roles of $f$ and $g$ interchanged. This completes the proof. 
        \vskip 5mm
        
         This work was supported by the Research Program P1-0291 from ARRS, Republic of Slovenia.
	        \vskip 10mm
	        \centerline{\bf REFERENCES}
	        \vskip 4mm
	        \rm
	        \noindent [AL]\ A.\ Alarc\'{o}n and F.\ J.\ L\'{o}pez:\ Null curves 
	        in $\C^3$ and Calabi-Yau conjectures.
	        
	        \noindent Math\ Ann.\ 355 (2013) 429-455
	        
	        \vskip 1mm
	        \noindent 
	        [AF]\ A.\ Alarc\'{o}n and F.\ Forstneri\v c:\ Every bordered Riemann surface 
	        is a complete proper curve in a ball.
	        
	        \noindent Math.\ Ann.\ 357 (2013) 1049-1070
	        \vskip 1mm
	        \noindent
	        [B]\ J.\ Bourgain: On the radial variation of bounded analytic functions on 
	        the disc.
	        
	        \noindent Duke Math.\ J.\ 69 (1993)671-682
	        \vskip 1mm
	        \noindent [J]\ P.\ W.\ Jones:\ A complete bounded complex submanifold of $\C^3$.
	        
	        \noindent Proc.\ Amer.\ Math.\ Soc.\ 76 (1979) 305-306
	        \vskip 1mm
	        \noindent [R]\ W.\ Rudin:\ \it Real and complex analysis,\rm \ third edition.   
	        
	        \noindent WCB/Mc Graw-Hill, New York, 1987
	        \vskip 1mm
	        \noindent [Y]\ P.\ Yang:\ Curvature of complex submanifolds of $\C^n$. 
	        
	        \noindent In: Proc.\ Symp.\ Pure.\ Math.\ Vol.\ 30, part 2, pp. 135-137. Amer.\ Math.\ Soc., 
	        Providence, R.\ I.\ 1977
	        \vskip 10mm
	        \noindent Institute of Mathematics, Physics and Mechanics
	        
	        \noindent Ljubljana, Slovenia
	        
	        \noindent josip.globevnik@fmf.uni-lj.si
	        
        \end